\documentclass[10pt,twocolumn,twoside]{article}
\usepackage[top=1in,bottom=1in,left=0.75in,right=0.75in,paperwidth=8.5in,paperheight=11in]{geometry}
\usepackage{fancyhdr}
\usepackage[utf8]{inputenc}
\usepackage{amsmath}
\usepackage{graphicx}
\usepackage{newtxmath}
\usepackage{authblk}
\usepackage[square,longnamesfirst,numbers]{natbib}

\newcommand{\diff}[1]{\mathrm{d}#1}
\newcommand{\diffop}{\mathrm{d}}

\newcommand{\dx}{\diff{x}}

\newcommand{\dy}{\diff{y}}

\title{Hyperreal Numbers for Infinite Divergent Series}
\author[1]{Jonathan Bartlett}
\author[2]{Logan Gaastra}
\author[3]{David Nemati}
\affil[1]{The Blyth Institute, jonathan.bartlett@blythinstitute.org}
\affil[2]{University of Michigan, logangaastra@gmail.com}
\affil[3]{The Blyth Institute}
\date{\today}

\begin{document}

\maketitle

\begin{abstract}
Treating divergent series properly has been an ongoing issue in mathematics.
However, many of the problems in divergent series stem from the fact that divergent series were discovered prior to having a number system which could handle them.
The infinities that resulted from divergent series led to contradictions within the real number system, but these contradictions are largely alleviated with the hyperreal number system.
Hyperreal numbers provide a framework for dealing with divergent series in a more comprehensive and tractable way.
\end{abstract}

\section{The Problem of Infinite Series}

Historically, infinities have led to many problems in mathematics.  
Infinities, when not handled carefully, easily lead to contradictions and indeterminacies.
Therefore, caution has always been urged when dealing with infinite series.

This is especially true with divergent infinite series.
Convergent infinite series generally behave unproblematically similar to the value that they converge to.
Given a series that converges to $2$ and another series that converges to $3$ then the sum of the values of the series will be $5$ and their product will be $6$.
Therefore, the nature of these series can be summarized into a single number.

With divergent series, this is not so straightforward.
A lack of agreement on the rules for handling infinities had led to numerous problems with handling divergent series.
If a series diverges to infinity, is it greater than or equal to some other series that diverges to infinity?
Can the terms of the series be rearranged?
Can their spacing be modified?
Is $1 + 1 + 1 + \ldots$ equivalent to $1 + 0 + 1 + 0 + 1 + 0 + \ldots$?

Lack of answers to questions like this have stifled work in divergent series, and have caused many mathematicians to think of divergent series as invalid entities to work with rigorously.

\section{Working with Infinities}

Many paradoxes exist with infinities.
For instance, are there the same number of positive even integers as positive integers?
There are an infinity of them, but does that make them the same?
It seems pretty obvious that, on a number line, positive integers occur twice as often.
However, there are an infinite amount of both.

Cantor's solution to this problem is to separate out the final quantity of a set (the cardinality) from the arrangment of a set (its ordinality).  
The cardinal numbers do not behave in any way similar to real numbers.
The ordinals, on the other hand, behave in many ways similar to real numbers.
However, Cantor's own system for ordinal arithmetic is difficult to use, and doesn't translate well between transfinite and regular real arithmetic.

The hyperreal number line has many similarities to Cantor's ordinals, operating essentially at the level of ``ordinal'' in Cantor's system.
However, the hyperreal number line offers a way to do arithmetic with infinities in a way that very closely matches real arithmetic through the use of the transfer principle \citep{henle2003}.
The transfer principle states that any first-order proposition that is true for the reals is also true for the hyperreals.
This means that the standard arithmetic principles for dealing with real numbers will apply to hyperreal numbers as well.

The hyperreal number line operates with an infinite unit, $\omega$, that represents an order of infinity.\footnote{The choice of character/typography for the unit varies with the author.  For instance, Keisler uses $H$ \citep{keisler2012}.  $\omega$ was chosen because of its historical connection with ordinal-type infinities.}
The way it is usually handled, $\omega$ isn't a specific number in the typical sense, but rather more of a benchmark of infinity.

Previous work has shown that hyperreal numbers could be a potential solution to how values of divergent series can be represented \citep{gaastra2016}.\footnote{Other work worth mentioning in this area are \citep{paterson2018a} and \citep{paterson2018b}.  In the current work, we will use a notation similar to \citep{keisler2012} to notate hyperreal values, and show how infinite series can be simplified to them.  Paterson did the opposite, by notating hyperreal values with the infinite sum that represents them.}
The present paper will build on this original idea and establish a system for using hyperreal numbers to assign values to infinite series.

\section{Hyperreals and Partial Sums}
\label{partialsum}

The vast majority of issues with divergent series comes with the transition from partial sums to infinity.
As long as a series remains a \emph{partial} sum, arithmetic with the series is unproblematic.
Therefore, it would be beneficial to develop a system which matched the partial sum behavior of finite sums, but allowed the result to be generalized to infinity.

The value of a partial sum of a given length is sensitive to the order of the terms in the infinite sequence.
Imagine summing the first $n$ terms of an infinite series.
The result will not be the same with different orderings of terms.
If the extent of the partial summation is unknown, then it is also unknown the extent to which numbers can be reordered.

Additionally, tacking on zeroes to the beginning of the series can potentially \emph{change} the partial sum.
Therefore, although adding zeroes to the beginning of a series has the \emph{appearance} of being a null operation, because doing so modifies the value of finite partial sums, it can also lead to long-term changes in behavior.

Therefore, just as ordinal infinities differ because of order-dependent properties, so too will infinite series exist as heavily order-dependent entities.

To understand many of the rules that will be developed for infinite series, imagine that the rules are being built for merely doing partial sums to an unknown parameter $k$, where $k$ at least acts like a particular finite value, but is larger than any particular list index referenced by any finite manipulation of the series.

Some of these formulas will be further reducible due to the nature of the hyperreals, as will be discussed in Section~\ref{secPrincipalValue}.

\section{Pinning Down $\omega$}

Since $\omega$ operates as a benchmark instead of a number, the first task is to identify the benchmark to associate $\omega$ with.
This is actually to some extent an arbitrary decision.
Any infinitely large value could be used to establish a baseline $\omega$.

However, the value that seems most natural for $\omega$ (especially for summation) is the size of the set of positive integers.
Therefore, $\omega$ will be used to refer to the total quantity of positive integers.

Because of this, the notation used will be more specific when writing summations.
Instead of summing to the ambiguous infinity, $\infty$, a summation to the specific infinity of all positive integers, $\omega$, will be used.
Therefore, the series $1 + 2 + 3 + \ldots$ will be written as
\begin{equation}
\sum_{i = 1}^{\omega} i
\end{equation}

This will establish the starting benchmark for relationships among the different series.

\section{The Standard Summation}

Because ordinal infinities are so order-dependent, it is important to establish an official standardization of summation.
That is, $\sum_{1}^{\omega}$ will be different from $\sum_{0}^{\omega}$.
Even though it looks like series with these types of sums will have an identical number of terms (after all they both have infinite terms), using this methodology the latter one will actually have more elements than the former.

This is due to the principle established in Section~\ref{partialsum}.  
If, instead of $\omega$ being infinite, pretend that $\omega$ was just an ordinary finite integer parameter.

Examine the series
\begin{equation}
\label{basicomega1}
\sum_{i = 1}^{\omega} 1.
\end{equation}
If $\omega$ represented an integer (say, $5$) instead of $\infty$, it would be obvious that this sum represents a different value from the series
\begin{equation}
\label{basicomega0}
\sum_{i = 0}^{\omega} 1.
\end{equation}
Equation~\ref{basicomega1} would represent the value $5$ while Equation~\ref{basicomega0} would represent the value $6$.
Therefore, it is clear that having matching indices matters.

In fact, our ability to sum divergent series will sometimes depend on having summations with equivalent numbers of terms.
Therefore, a ``standard'' starting point for summation will need to be established in order to ensure that like entities are being compared and reasoned about.
In computations it technically winds up not mattering whether the starting point is $1$ or $0$, though the formulas would have to be reworked based on the starting index. 
However, since $\omega$ has been defined as being the size of the set of all positive integers, it makes sense to start at $1$.
For the purposes of this paper, the ``standard'' way of summing will be to start with $1$ and proceed to $\omega$.

\section{Simple Arithmetic and Geometric Series}
\label{secSimpleSeries}

\subsection{Arithmetic Series}

Arithmetic series take the form
\begin{equation}
\sum_{i = 1}^{n} a + (i - 1)d.
\end{equation}
The sum of an arithmetic series, given a starting value $a$, the number of elements $n$, and distance between elements $d$, can be given by the formula
\begin{equation}
\sum_{i=1}^{n} a + (i - 1)d = \frac{n}{2}\left(2a + (n - 1)d\right).
\end{equation}
To find the sum of an infinite arithmetic series, $\omega$ is used for $n$, forming a hyperreal value.
That reduces the formula to
\begin{equation}
\label{arithmeticSeriesFormula}
\sum_{i=1}^{\omega} a + (i - 1)d = \omega a + \frac{\omega^2d}{2} - \frac{\omega d}{2}.
\end{equation}

Therefore, to find the summation of the series $1 + 1 + 1 + \ldots$, one must only substitute in the correct parameters.
Since the starting value is $1$ and the distance between terms is $0$, this yields
\begin{align}
\sum_{i=1}^{\omega} 1 &= \omega \cdot 1 + \frac{\omega^2 \cdot 0}{2} - \frac{\omega \cdot 0}{2} \\
	&= \omega + 0 - 0\\
	&= \omega.
\end{align}
It is intuitively obvious that since there are $\omega$ $1$s added together that the sum of them would add up to $\omega$, as would be true for any finite value as well.

The arithmetic series $1 + 2 + 3 + \ldots$ can be calculated using hyperreals as well.
\begin{align}
\sum_{i=1}^{\omega} i &= \omega \cdot 1 + \frac{\omega^2 \cdot 1}{2} - \frac{\omega \cdot 1}{2} \\
 &= \frac{\omega^2}{2} + \frac{\omega}{2}.
\end{align}

The next arithmetic series to examine is $1 + 3 + 5 + \ldots$, which can be similarly calculated.
\begin{align}
\sum_{i=1}^{\omega} (2i - 1) &= \omega \cdot 1 + \frac{\omega^2 \cdot 2}{2} - \frac{\omega \cdot 2}{2} \\
 &= \omega^2.
\end{align}
Thus, the value of $1 + 3 + 5 + \ldots$ is equal to $(1 + 1 + 1 \ldots)^2$.

Interestingly, as noted in Section~\ref{partialsum}, there is nothing intrinsically infinite about the behavior of $\omega$ in these series.  
For instance, if $\omega$ was replaced with $5$, the results would hold.
That is, $(1 + 1 + 1 + 1 + 1)^2 = (1 + 3 + 5 + 7 + 9) = 25$.

Even though the sums are divergent, summing them has a very well-defined behavior within the combined hyper-real/partial sum methodology presented here.

\subsection{Geometric Series}

Geometric series take the form
\begin{equation}
\sum_{i = 1}^{n} a r^{i - 1},
\end{equation}
where $n$ is the number of terms, $a$ is the starting term, and $r$ is the common ratio.

A value for a geometric series can be given by the formula
\begin{equation}
\sum_{i = 1}^{n} a r^{i - 1} = a \frac{1 - r^n}{1 - r}.
\end{equation}
Because an infinite series will have $\omega$ terms, $n$ can be replaced with $\omega$.

Let us begin by looking at the series $1 + 2 + 4 + 8 + \ldots$.
The value of this series can be given by the formula
\begin{align}
\sum_{i = 1}^{\omega} 2^{i - 1} &= 1\cdot\frac{1 - 2^\omega}{1 - 2} \\
 &= 2^\omega - 1.
\end{align}
Divergent geometric series will generally have the same form.

Convergent series are also interesting.
The series $1 + \frac{1}{2} + \frac{1}{4} + \ldots$ can be plugged into the formula to yield
\begin{align}
\sum_{i = 1}^{\omega} \frac{1}{2}^{i - 1} &= 1\cdot\frac{1 - \frac{1}{2}^\omega}{1 - \frac{1}{2}} \\
 &= 2 - 2\cdot \left(\frac{1}{2}\right)^\omega \label{simpleConvergentGeometric}
\end{align}

\section{Generalizing to the Principal Value}
\label{secPrincipalValue}

In most discussions of hyperreal numbers, the halo of a number is considered the hyperreal values which are infinitely close to a standard real number.
However, this definition is too focused on real numbers.  

We will consider the order of a hyperreal value to be its largest exponent of $\omega$.
This is the most significant term of the hyperreal value.
We will call this most significant term the principal value of the hyperreal.
The halo (also known as a monad) of a hyperreal consists of all of the hyperreals which have the same principal value.\footnote{Most texts on hyperreal numbers define the halo or monad of $x$ to be all of the values $y$ for which $x - y$ is infinitesimal \citep[pg.~21]{loeb2015} \citep[pg.~52]{goldblatt1998}.  
However, defined in such a way, the infinitesimals $\omega^{-1}$ and $2\omega^{-1}$ are within a monad.  
Using principal values, $\omega^{-1}$ and $2\omega^{-1}$ are in the same galaxy, but not the same monad.  
You would have to have a term of lower-order infinity to be within a monad, such as $\omega^{-1}$ and $\omega^{-1} + \omega^{-2}$.  
This seems to be the essence of what the other texts are getting at, but, since most mathematics focuses on the reals, their definitions were entirely based on using reals as a starting point.  
Here, since we will have results in the hyperreals, we need definitions that are equally useful when the final result is a hyperreal number.}
We will use the $\simeq$ operator to denote two hyperreals which share the same principal value.\footnote{In practice, $\simeq$ can be replaced with $=$, as it denotes equality to the extent normally practiced in mathematics.  For instance, the differential $\diffop\left(xy\right)$ is often stated as being \emph{equal} to $x\,\dy + y\,\dx$, but really it is just the principal value.  The actual value is $x\,\dy + y\,\dx + \dy\,\dx$.  The $\dy\,\dx$ term is always discarded because it is infinitely less significant than the other pieces.  Even when discarding this term, the equality sign is used.  Therefore, while the present paper will be pedantic about asserting exact equality or mere principal value, for most general purposes equality can be asserted even when only stating the principal value.}
Therefore, the halo of a hyperreal number consists of all of those numbers which share the same principal value.

Many people use ``infinitely close'' as a colloquialism to describe two hyperreals which share the same principal value.
However, technically it is not correct, since, when dealing with infinities, two hyperreals which differ by multiple infinities can be considered ``infinitely close.''
That is, if $d \neq 0$, then $\omega^2 + 5\omega$, $\omega^2 - 12\omega$, and $\omega^2 + 23$ all share the same principal value, $\omega^2$.  
They are infinitely apart, yet, colloquially, they can be considered ``infinitely close'' because their differences are infinitely less significant than their similarities.

When dealing with hyperreals, the principal value is the main one of concern.
So, for instance, while $1 + 2 + 3 + \ldots$ is exactly described by $\frac{\omega^2}{2} + \frac{\omega}{2}$, its principal value is just $\frac{\omega^2}{2}$.
Therefore, the formula given in (\ref{arithmeticSeriesFormula}) can actually be simplified to
\begin{equation}
\sum_{i=1}^{\omega} a + (i - 1)d \simeq \frac{\omega^2d}{2}
\end{equation}
if $d \neq 0$.\footnote{When $d = 0$, then the $\omega^2$ term goes to zero, and the series simplifies to $a\cdot \omega$ instead.}

Interestingly, we can see that, while the exact value of the hyperreal associated with a series depends on the starting point, the principal value depends only on the distance chosen, provided that $d \neq 0$.

Geometric series can use similar considerations.
You may have noticed that the hyperreal given for the series $1 + \frac{1}{2} + \frac{1}{4} + \ldots$ in Section~\ref{secSimpleSeries} is $2 - 2\cdot \left(\frac{1}{2}\right)^\omega$.
Typically, this series is thought to converge to $2$.
In fact, its principal value is $2$, because $\left(\frac{1}{2}\right)^\omega$ is an infinitesimal.

The use of principal values allows for a great amount of simplification for hyperreal values and formulas.

As an example, the ratio between two given arithmetic series can be solved for very simply.
\begin{align}
S_1 &= \sum_{i=1}^{\omega} a_1 + (i - 1)d_1 \simeq \frac{\omega^2(d_1)}{2} \nonumber \\
S_2 &= \sum_{i=1}^{\omega} a_2 + (i - 1)d_2 \simeq \frac{\omega^2(d_2)}{2} \nonumber \\
\frac{S_1}{S_2} &\simeq \frac{\frac{\omega^2(d_1)}{2}}{\frac{\omega^2(d_2)}{2}} = \frac{d_1}{d_2}
\end{align}
In other words, the principal value of the ratio of two arithmetic series is simply the ratio of the distances.

\section{Series Manipulation Rules for Finite Subsets}


Many attempts to manipulate divergent series have resulted in contradictions, to the extent that many suggest that it is best to not attempt to do so.
The reason for these contradictions, however, lies in the treatment of the infinite nature of the number of values.

In the real system, $\infty$ is considered a boundless number.
That is, there is not $\infty + 1$ that is distinct from $\infty$.
Likewise, $\infty - 1$ is also infinity.
Essentially, within the real numbers, $\infty$ is used largely like an ambiguous infinite value, essentially saying that ``the real numbers can't handle this value.''

If, instead, the hyperreal numbers are used, then $\omega$ and $\omega + 1$ are distinct quantities, despite the fact that they are both infinite.
The rules for manipulating series come from these ideas.

\subsection{Finite Term Addition}
\label{finitetermaddition}

To begin with, it is possible to easily add a scalar value to a series, provided that it is added \emph{to} one of the particular terms of the series.
In other words, suppose the value $A$ is added to the series $1 + 2 + 3 + \ldots$.
This can be written as
\begin{equation}
A + \sum_{i = 1}^{\omega} i
\end{equation}
or as 
\begin{equation}
A + (1 + 2 + 3 + \ldots).
\end{equation}
To integrate $A$ into the series, $A$ can be added to any distinct position.
The series could read as 
\begin{equation}
(A + 1) + 2 + 3 + \ldots
\end{equation}
or
\begin{equation}
1 + 2 + (A + 3) + \ldots.
\end{equation}
All of these yield the same value for the final series, as long as partial sums are taken starting after the index where $A$ is added.

Additionally, $A$ can be spread across multiple finite terms.
For instance, half of $A$ could be added to each of the first two terms, yielding
\begin{equation}
A + (1 + 2 + 3 + \ldots) = \left(1 + \frac{A}{2}\right) + \left(2 + \frac{A}{2}\right) + 3 + \ldots .
\end{equation}
In fact, there is no reason why the same amount would have to be distributed to each position.
\begin{equation}
A + (1 + 2 + 3 + \ldots) = \left(1 + \frac{2}{5}A\right) + \left(2 + \frac{3}{5}A\right) + 3 + \ldots
\end{equation}

\subsection{Finite Term Insertion and Removal}
\label{terminsertion}

Because this method of summation is based on partial sums, it should be apparent that inserting and removing terms will in fact alter the summation.
For instance, let's begin with the arithmetic sum $1 + 1 + 1 + \ldots$.
It may seem intuitive that one should be able to freely add or remove a $1$ from this sum without affecting the sum.
In this particular series, the exact hyperreal value does change, but not the principal value.

Again, remember that, as mentioned in Section~\ref{partialsum}, this conception of summation will be based on partial sums.
So, let us begin by considering the partial sum
\begin{equation}
\sum_{i = 1}^{k} 1.
\end{equation}
If $k$ is a finite number, then adding one to this sequence will in fact alter its value.
Additionally, removing a $1$ from this sequence will also alter its value.
Therefore, 
\begin{equation}
\sum_{i = 1}^{k} 1 \neq 1 + \sum_{i = 1}^{k} 1.
\end{equation}
Likewise,
\begin{equation}
\sum_{i = 1}^{k} 1 \neq \sum_{i = 0}^{k} 1 \neq \sum_{i = 2}^{k} 1.
\end{equation}
Because performing these operations will change the value for any partial sum of $k$ terms for a finite $k$, they will also change the value for a hyperreal $k$ such as $\omega$.
However, for these particular series, the principal value will be the same, because $\omega \simeq \omega + 1 \simeq \omega - 1$.

Additionally, a more surpising fact is that removing a term from a sequence also changes its value if it does not also change the number of terms being summed.
Consider the series
\begin{equation}
\label{seriesunbumped}
1 + 2 + 3 + \ldots = \sum_{i = 1}^{\omega} i.
\end{equation}
This series is not equal to the series
\begin{equation}
\label{seriesbumped}
1 + \sum_{i = 1}^{\omega} (i + 1).
\end{equation}
although it does have the same principal value in this case.
In other words,
\begin{equation}
(1 + 2 + 3 + \ldots) \neq 1 + (2 + 3 + 4 + \ldots)
\end{equation}
but
\begin{equation}
(1 + 2 + 3 + \ldots) \simeq 1 + (2 + 3 + 4 + \ldots).
\end{equation}

The reason for this is readily apparent when considering how these work in terms of partial sums.
If the parameter $k$ was used instead of $\omega$, then it is apparent that the value of (\ref{seriesbumped}) actually has an \emph{extra} term compared to (\ref{seriesunbumped}).
That is, it is obvious that
\begin{equation}
\sum_{i = 1}^{5} i \neq 1 + \sum_{i = 1}^5 (i + 1).
\end{equation}
This can also be seen in the results of applying the arithmetic series formula to the two series.
For $(1 + 2 + 3 + \ldots)$ the formula yields $\frac{\omega^2}{2} + \frac{\omega}{2}$.
However, for $(2 + 3 + 4 + \ldots)$ the formula yields $\frac{\omega^2}{2} + \frac{3}{2}\omega$.

Now, terms can be removed without even affecting the exact hyperreal value if they are replaced by zeroes in the sequence, or if the sequence starting index is moved appropriately.
In other words,
\begin{equation}
\label{termremoval}
(1 + 2 + 3 + \ldots) = 1 + (0 + 2 + 3 + \ldots) = 1 + \sum_{i = 2}^{\omega} i.
\end{equation}
This can be easily proved using the principle derived in Section~\ref{finitetermaddition}.
For instance, to move the $1$ outside of the series, $1 + -1$ can be added to the series.
\begin{align}
1 + -1 + (1 + 2 + 3 + \ldots)  \nonumber \\
  ~~~~~~ &= 1 + ( (1 + -1) + 2 + 3 + \ldots)  \nonumber \\
  &= 1 + (0 + 2 + 3 + \ldots)
\end{align}

\subsection{Finite Term Rearrangement}

As can be deduced from Sections~\ref{finitetermaddition} and \ref{terminsertion}, any number of finite terms in a series can be rearranged in position.
That is, for any given series member with a value of $A$, $A - A$ can be added to the series, applying the $-A$ such that it cancels out the value of the series member.
After doing this to several series members, the inverse operations can then be applied to move these values to any finite position in the series.

Doing this will preserve the partial summing behavior of the series for all partial sums after the members which have been manipulated.

\section{More Advanced Series}
\label{advseries}

While basic formulas for divergent series of arithmetic and geometric series can be established using the standard formulas, more advanced series require the use of discrete integral calculus to establish the formulas for series.
Doing so leads to very interesting results.

\subsection{Ces\`aro Sums and Oscillating Series}

Oscillating series have an interesting history of treatment within mathematics.
The standard series to consider is Grandi's series: $1 - 1 + 1 - 1 + \ldots$.  
Or, written more formally,
\begin{equation}
\sum_{i = 1}^{\infty} (-1)^{i + 1}.
\end{equation}
Partials sums for this series can be found by performing a discrete integral.
\begin{equation}
\sum_{i = 1}^{n} (-1)^{i + 1} = \frac{1}{2}(-1)^{n + 1} + \frac{1}{2}.
\end{equation}
What is particularly interesting about this formula is that the Ces\`aro sum  of the infinite series ($\frac{1}{2}$) is present in the formula.

Now, consider the oscillating series $-1 + 1 - 1 + \ldots$.
This series has the formula
\begin{equation}
\sum_{i = 1}^{\infty} (-1)^i.
\end{equation}
A discrete integral of the partial sums yields the formula
\begin{equation}
\sum_{i = 1}^{n} (-1)^i = \frac{1}{2}(-1)^n - \frac{1}{2}.
\end{equation}
Note that in this as well, $-\frac{1}{2}$ is the Ces\'aro summation of the infinite series.

This leads to the conjecture that, in evaluating infinite series using integral formulas, 
\begin{equation}
\label{negoneconjecture}
(-1)^\infty = 0,
\end{equation}
at least for additive offsets of $\omega$.
For instance, in the case of Grandi's series, using the $\omega$ notation, the infinite series would include $(-1)^{\omega + 1}$.  
The other series includes $(-1)^{\omega}$.
According to the present conjecture, both of these simplify to $0$, at least for the purpose of creating formulas for infinite series based on partial sums.

This can be understood probabilistically.
Since we have no information about what sign $-1^{\omega}$ will have, we can say that 
\begin{equation}
    -1^{\omega} = \pm 1.
\end{equation}
Since both of these possibilities are equally probable, the limit towards infinity resolves to their average, or zero.
Also, since we have no information about the sign of $-1^{\omega}$, we have equally little information about the sign of $-1^{\omega + 1}$, or any other variation on $\omega$ which is not biased towards evenness (e.g., $2\omega$).

The expression $-1^{x}$ has an oscillation pattern very similar to $\sin(x)$.
Since \citep{paterson2018a} showed that $\sin(\omega) = 0$ in the surreal numbers, it is possible that a similar proof may be found for $-1^\omega = 0$ along similar lines in the hyperreals.

\subsection{Other Oscillatory Behavior}

Because (a) discrete integration can be used to find formulas for series involving partial sums, and (b) the behavior of $(-1)^\infty$ (for infinities without bias towards evenness) is conjectured to be zero, the behavior of a wide variety of oscillatory behaviors can be deduced.

Raising $-1$ to the $i$th power can produce all sorts of oscillatory behavior.
As has been seen with Grandi's series, this can produce a series of values that go back-and-forth across a mean value (the mean value can be changed by adding, and the back-and-forth can be changed by multiplying).

However, $(-1)^i$ can also be expanded to blank out members of a series.
For instance, to blank out every other member of a series, the formula
\begin{equation}
\label{blankingformula}
\frac{((-1)^i + 1)}{2}
\end{equation}
can be used.
This simplifies to $1$ where $i$ is even and $0$ when $i$ is odd.
Therefore, by multiplying a given formula by (\ref{blankingformula}), pieces of the given formula will be zeroed out.

For instance, take the series $1 + 2 + 3 + \ldots$.
This series can be converted to the series $0 + 2 + 0 + 4 + 0 + 6 + \ldots$ by applying (\ref{blankingformula}).
This gives the series
\begin{equation}
\sum_{i = 1}^{\omega} i\cdot\left(   \frac{((-1)^i + 1)}{2} \right).
\end{equation}
The discrete integral yields
\begin{equation}
\label{otheroscillator}
\sum_{i = 1}^{n} i\cdot\left(\frac{((-1)^i + 1)}{2}\right) = \frac{1}{8}\left(2n^2 + 2n(-1)^n + 2n + (-1)^n - 1\right)
\end{equation}
When $n = \omega$ the formula runs into a problem with simplifying this through the conjecture (\ref{negoneconjecture}) because it yields an indeterminate form.
The term $2n(-1)^n$ becomes an indeterminate form of the type $\omega\cdot 0$.
This can be resolved, however, through L'Hospital's Rule.
\begin{equation}
\label{lhospital}
\lim_{n\to\infty} \frac{2n}{(-1)^{-n}} = \frac{2}{-\ln(-1)(-1)^{-n}} = -\frac{2}{\ln(-1)}(-1)^n.
\end{equation}
Now (\ref{negoneconjecture}) can be applied without ambiguity, simplifying it to zero.

Therefore, for $n = \omega$, (\ref{otheroscillator}) simplifies to
\begin{equation}
\sum_{i = 1}^{n} i\cdot\left(\frac{((-1)^i + 1)}{2}\right) = \frac{1}{8}\left(2n^2 + 2n - 1\right).
\end{equation}

This means that the value of this sum in the hyperreals is $\frac{1}{4}\omega^2 + \frac{1}{4}\omega - \frac{1}{8} \simeq \frac{1}{4}\omega^2$.

Interestingly, this is a different result than for the simple series $2 + 4 + 6 + \ldots$.
Since $2 + 4 + 6 + \ldots$ is a simple arithmetic series, we can determine the hyperreal sum using (\ref{arithmeticSeriesFormula}).
\begin{equation}
\sum_{i = 1}^{\omega} 2 + (i - 1)2 = \omega^2 + \omega \simeq \omega^2.
\end{equation}
This is a different result than what was obtained for $0 + 2 + 0 + 4 + 0 + 6 + \ldots$, which was $\frac{1}{4}\omega^2$, indicating that the two series have different behaviors.

\subsection{$1 - 2 + 3 - 4 + \ldots$}

Euler's sum for the series $1 - 2 + 3 - 4 + \ldots$ can be confirmed using this method as well.
This series can be given the value
\begin{equation}
\sum_{i = 1}^n i(-1)^{i - 1} = \frac{1}{4}\left(-2n(-1)^n + (-1)^{n + 1} + 1\right).
\end{equation}
Using (\ref{negoneconjecture}) and (\ref{lhospital}) this simplifies to $\frac{1}{4}$.

Interestingly, this is one of the few functions that is not changed even in its exact hyperreal by prepending a zero to the function.
\begin{equation}
\sum_{i = 1}^n (i - 1)(-1)^{i} = \frac{1}{4}\left(2n(-1)^n + (-1)^{n + 1} + 1\right).
\end{equation}
Likewise, (\ref{negoneconjecture}) allows this to reduce to $\frac{1}{4}$.

\section{Whole Series Manipulation Rules}

In addition to manipulation of finite partial sums of a series, certain operations can (and can't) be performed to the series as a whole.
In this section, some of these operations will be considered.

\subsection{Scalar Multiplication}

Because of the distributivity of multiplication, multiplication of a series by a scalar value will distribute the scalar multiplication to every term.
\begin{equation}
2(1 + 2 + 3 + \ldots) = (2\cdot 1 + 2\cdot 2 + 2\cdot 3 + \ldots).
\end{equation}
Or, written as a formula,
\begin{equation}
n \sum_{i = 1}^{\omega} f(i) = \sum_{i = 1}^{\omega} n\,f(i).
\end{equation}

\subsection{Whole Series Addition}
\label{seriesaddition}

Adding two series together is equivalent to a term-by-term addition of the series.
Since the method presented here is based on partial sums, term-by-term addition only works when the lower and upper bounds of the terms are identical.

Therefore,
\begin{equation}
\left( \sum_{i = 1}^{\omega} f(i) \right) + \left( \sum_{i = 1}^{\omega} g(i) \right) = \sum_{i = 1}^{\omega} f(i) + g(i).
\end{equation}
However,
\begin{equation}
\label{additionlimitsneq}
\left( \sum_{i = 0}^{\omega} f(i) \right) + \left( \sum_{i = 1}^{\omega} g(i) \right) \neq \sum_{i = 1}^{\omega} f(i) + g(i)
\end{equation}
because the limits of summation differ.
Again, to see why this is the case, imagine replacing $\omega$ with a fixed scalar such as $5$.
In (\ref{additionlimitsneq}), the left-hand addend would have a different number of terms than the right-hand addend.

\subsection{Series Spacing}

As noted in Section~\ref{terminsertion}, adding or removing elements of a series, even if they are zero, has an effect on the sum of the resulting series.
This effect can be calculated using the considerations discussed in Section~\ref{advseries}.

For instance, the series $1 + 1 + 1 + \ldots$ can be spaced out by adding in zeroes, to make $1 + 0 + 1 + 0 + \ldots$.
A variation of the oscillatory pattern in (\ref{blankingformula}) can be used to give the series the formula
\begin{equation}
\sum_{i = 1}^{n}  \frac{((-1)^{i + 1} + 1)}{2}.
\end{equation}

The discrete integral of this yields the formula
\begin{equation}
\label{seriesoneblanksright}
\frac{1}{2}n + \frac{1}{4}(-1)^{n + 1} + \frac{1}{4}
\end{equation}
Using conjecture (\ref{negoneconjecture}) this reduces to the hyperreal value $\frac{1}{2}\omega + \frac{1}{4} \simeq \frac{1}{2}\omega$.

This is a slightly different value (but with the same principal value) than for the series $0 + 1 + 0 + 1 + \ldots$.
This series can be represented as
\begin{equation}
\label{seriesoneblanksleft}
\sum_{i = 1}^{n}  \frac{((-1)^{i} + 1)}{2} = \frac{1}{2}n + \frac{1}{2}(-1)^n - \frac{1}{4}.
\end{equation}
Using conjecture (\ref{negoneconjecture}), the hyperreal value for this is $\frac{1}{2}\omega - \frac{1}{4} \simeq \frac{1}{2}\omega$.

If (\ref{seriesoneblanksright}) and (\ref{seriesoneblanksleft}) were added, it should be equivalent whether they are added term-by-term (Section~\ref{seriesaddition}) or by summing their relevant values.

Summing term-by-term it is apparent that
\begin{align}
&(1 + 0 + 1 + 0 + \ldots) \nonumber \\
+ ~~ &(0 + 1 + 0 + 1 + \ldots) \nonumber \\
= ~~ &(1 + 1 + 1 + 1 + \ldots).
\end{align}
The value of this series was deduced to be $\omega$ in (\ref{basicomega1}).
Likewise, if the values for each series are added the result is
\begin{equation}
\left(\frac{1}{2}\omega + \frac{1}{4}\right) + \left(\frac{1}{2}\omega - \frac{1}{4}\right) = \omega.
\end{equation}

\section{Conclusion}

Here a method of summation was presented that uses the structure of the hyperreal numbers to represent values for divergent series.
This methodology was shown to be stable across a variety of different scenarios.
One unproven, but seemingly correct, conjecture was relied upon for this formulation.
Future work will focus on proving (\ref{negoneconjecture}).

\section{Acknowledgements}

I wanted to take a moment to thank Stanley Schmidt.
I was thinking on this problem at the same time I was reading his \textit{Life of Fred} books to my children.
The fundamental idea for this method of summation came from thinking about Gaastra's original presentation \citep{gaastra2016} while reading \textit{Life of Fred: Kidneys} to my children, when Fred was using the formula for arithmetic series \citep{fredkidneys}.

Additionally, \textit{The Infinite} by A. W. Moore provided some help to the imagination in his discussion of the L\"owenheim-Skolem theorem.
The basics of the discussion was to point out that there was little in the theory of infinities that were really unique to infinity.
Even finite sets can look ``infinite'' in some ways to other sets.
The techniques and ideas explored in Section~\ref{partialsum} were based largely off of thinking about infinities as much more tame and finite-like than is normally considered.

Finally, I want to thank Jessica Hastings, whose interest in the ``Wheat and Chessboard'' problem\citep{weissteinwheat} originally introduced me to the concepts in discrete calculus.

\bibliographystyle{ieeetr}
\bibliography{DivergentSeriesHyperreals}
\end{document}